\documentclass[12pt,reqno]{amsart}

\usepackage[dvips]{graphicx}
\usepackage[arrow,matrix,curve]{xy}

\usepackage{amssymb, latexsym, amsmath, amscd, array, 
}

\newcommand\N {{\mathbb N}} 
\newcommand\R {{\mathbb R}}
\newcommand\Q{{\mathbb Q}}

\newcommand\RRR{{\mbox{I\!I\!R}}}

\newcommand\Los{{\L}o{\'s}}
\newcommand\fnull{\mathcal{F}_{\!n\!u\!l\!l}}
\newcommand\fez{\mathcal{F}_{ez}}

\theoremstyle{definition}

\numberwithin{equation}{section}
\numberwithin{figure}{section} 
\numberwithin{table}{section}

\begin{document}

\title{A Cauchy--Dirac delta function}

\author[M.~Katz]{Mikhail G. Katz} \address{Department of Mathematics,
Bar Ilan University, Ramat Gan 52900 Israel}
\email{katzmik@macs.biu.ac.il}

\author{David Tall} \address{Mathematics Education Research Centre,
University of Warwick, CV4 7AL, United Kingdom}
\email{david.tall@warwick.ac.uk}



\date{\today}

\begin{abstract}
The Dirac $\delta$ function has solid roots in 19th century work in
Fourier analysis and singular integrals by Cauchy and others,
anticipating Dirac's discovery by over a century, and illuminating the
nature of Cauchy's infinitesimals and his infinitesimal definition of
$\delta$.
\end{abstract}

\maketitle

\tableofcontents

\section{Dirac on infinity}
\label{delta}

The specialisation of the scientific disciplines since the 19th
century has led to a schism between the scientists' pragmatic
approaches (e.g., using infinitesimal arguments), on the one hand, and
mathematicians' foundational preferences, on the other.  The widening
schism between mathematics and physics has been the subject of much
soul-searching, see, e.g., S.P. Novikov~\cite{No2}.

In an era prior to such a schism, a seminal interaction between
physics and mathematics was foreshadowed in a remarkable fashion,
exploiting infinitesimals, in Cauchy's texts from 1827.  Cauchy
defined the so-called Dirac delta function in terms of infinitesimals,
and applied it to the evaluation of singular integrals and in Fourier
analysis.

Let us first review P.~Dirac's well-known discussion of~$\delta$.
Dirac writes in his \textsection~10:
\begin{quote}
Consider now the single integral
\[
\hskip1.3in \int \langle \xi'x\vert\xi''y \rangle d\xi''.\hskip1.3in
(29)
\]
From the orthogonality theorem, the integrand here must vanish over
the whole range of integration except the one point~$\xi''=\xi'$.
[\ldots] in general~$\langle \xi'x\vert\xi'y \rangle$ must be
infinitely great in such a way as to make (29) non-vanishing and
finite.  The form of infinity required for this will be discussed in
\textsection~15 (Dirac 1958, \cite[p.~39]{Dir}).
\end{quote}
Note that Dirac explicitly speaks of ``infinitely great'' values of
the integrand in his formula (29), and does not shy away from
discussing ``the form of infinity required''.  Dirac resumes his
discussion of~$\delta$ in his \textsection~15 in the following terms:
\begin{quote}
Our work in \textsection~10 led us to consider quantities involving a
certain kind of infinity.  To get a precise notation for dealing with
these infinities, we introduce a quantity~$\delta(x)$ depending on a
parameter~$x$ satisfying the conditions
\[
\begin{cases}
\int_{-\infty}^{\infty} \delta(x)dx=1 \\ \delta(x)=0 \; \text{for} \;
x\not=0
\end{cases}
\]
(Dirac 1958, \cite[p.~58]{Dir}).
\end{quote}
Dirac views his~$\delta$ as a way of ``dealing with these
infinities''.  Dirac is not content with merely providing such an
algebraic definition, and proceeds to explain how one can ``get a
picture of~$\delta(x)$'':
\begin{quote}
take a function of the real variable~$x$ which vanishes everywhere
except inside a small domain, of length~$\epsilon$ say, surrounding
the origin~$x=0$, and which is so large inside this domain that its
integral over this domain is unity.%
\footnote{\label{spec1}Dirac does not specify a function, merely
mentioning that the procedure works ``provided there are no
unnecessarily wild variations''.  A century earlier Cauchy did specify
a function (see formula~\eqref{31b} and footnote~\ref{f7}).}
[\ldots] Then in the limit~$\epsilon\to 0$ this function will go over
into~$\delta(x)$ (ibid.).
\end{quote}
Here Dirac is using the expression ``in the limit'' in a generic sense.%
\footnote{Dirac's ``limit'' can be interpreted mathematically either
in terms of distribution theory \`a la Sobolev and Schwartz, or in the
sense of an ultraproduct.}
Furthermore,
\begin{quote}
The most important property of~$\delta(x)$ is exemplified by the
following equation,
\[
\int_{-\infty}^{\infty} f(x)\delta(x)=f(0),
\]
where~$f(x)$ is any continuous function of~$x$ (Dirac
\cite[p.~59]{Dir}).
\end{quote}
The above equation is identical to one appearing in Cauchy a century
earlier.%
\footnote{\label{f3}See formula~\eqref{151}.}
An additional intriguing formula appears in the context of Dirac's
discussion of the discontinuity of the principal value of the log
function:
\begin{equation}
\label{11b}
\frac{d}{dx} \log x = \frac{1}{x} - i\pi \delta(x),
\end{equation}
accompanied by a comment to the effect that ``this particular
combination of reciprocal function and~$\delta$ function plays an
important role in the quantum theory of collision processes''%
\footnote{The imaginary part of the principal value of~$\log(x)$
in~\eqref{11b} is~$0$ for~$x>0$ and~$i\pi$ for~$x<0$.  Thus,
$\operatorname{Im}(\log(x))=-i\pi H(x)$ where $H(x)$ is the Heaviside
function (which equals~$-1$ for negative~$x$ and~$0$ for
positive~$x$), whose derivative is~$\delta(x)$ in the sense discussed
in Section~\ref{nine}.}
(Dirac \cite[p.~61]{Dir}).

\section{Modern Cauchy scholarship}

The earliest mention of Cauchy's~$\delta$ appears to be in
H.~Freudenthal.

\subsection{Freudenthal's biography}

Cauchy's anticipations of~$\delta$ were alluded to in 1971 by
Freudenthal, who mentioned~$\delta$ functions twice in his Cauchy
biography for the Dictionary of Scientific Biography \cite{Fr}.  They
were the subject of scholarly attention by J.~L\"utzen (1982,
\cite{Lut82}) and D.~Laugwitz~(1989, \cite{Lau89}; 1992 \cite{Lau92}).
A spate of recent articles have undertaken a re-evaluation of Cauchy's
foundational contribution to the development of infinitesimal
analysis: Br\aa ting (2007, \cite{Br}), Barany (2011, \cite{Ba11}),
B\l aszczyk et al. (2012, \cite{BKS}), Borovik et al. (2011,
\cite{BK}), Katz \& Katz~\cite{KK11a, KK11b}.  The present text will
focus on Cauchy's anticipations of the Dirac delta.

J.~L\"utzen \cite{Lut82} traces the origins of distribution theory
(traditionally attributed to S.~Sobolev and L.~Schwartz), in 19th
century work of Fourier, Kirchhoff, and Heaviside.  An accessible
introduction to this area may be found in J.~Dieudonn\'e~\cite{Die}.

\subsection{Dieudonn\'e's review}
\label{12}

Dieudonn\'e is one of the most influential mathematicians of the 20th
century.  A fascinating glimpse into his philosophy is provided by his
1984 review \cite{Die} of J.~L\"utzen's book.  At the outset,
Dieudonn\'e poses the key question:
\begin{quote}
One [\ldots] may well wonder why it took more than 30 years for
distribution theory to be born, after the theory of integration had
reached maturity (Dieudonn\'e \cite[p.~374]{Die}).
\end{quote}
This remark inscribes itself in a tradition of a long-established and
pervasive dogma, to the effect that 20th century physicists were the
first to invent the delta function.  Thus, M.~Bunge, responding to
Robinson's lecture {\em The metaphysics of the calculus\/}, evoked the
physicist's custom
\begin{quote}
to refer to the theory of distributions for the legalization of the
various delta `functions' which his physical intuition led him to
introduce (Bunge in Robinson, 1967 \cite[p.~44-45]{Ro67};
\cite[p.~553-554]{Ro79}).
\end{quote}
But did Dirac introduce the delta function?
Laugwitz~\cite[p.~219]{Lau89} notes that probably the first appearance
of the (Dirac) delta function is in the~1822 text by
Fourier~\cite{Fo}.

In his review of J.~L\"utzen's book for Math Reviews, F. Smithies
notes:
\begin{quote}
Chapter 4, on early uses of generalized functions, covers fundamental
solutions of partial differential equations, Hadamard's ``partie
finie'', and many early uses of the delta function and its
derivatives, including various attempts to create a rigorous theory
for them (Smithies, \cite{Sm}).
\end{quote}

At the end of his review, Smithies mentions Cauchy: ``In spite of the
thoroughness of his coverage, [L\"utzen] has missed one interesting
event---A. L. Cauchy's anticipation of Hadamard's `partie finie' '',
but says not a word about Cauchy's infinitesimally-defined delta
functions.

\subsection{Gilain on limits}

Cauchy's use of infinitesimals is not always reported accurately in
the literature.  Thus C.~Gilain claims that
\begin{quote}
On sait que Cauchy d\'efinissait le concept d'infiniment petit \`a
l'aide du concept de limite, qui avait le premier r\^ole%
\footnote{Translation: ``We know that Cauchy defined the concept of
infinitely small by means of the concept of limit, which played the
primary role'' (Gilain).}  (voir Analyse alg\'ebrique, p.~19 \dots)
(Gilain, 1989, \cite[footnote~67]{Gil}).
\end{quote}
Here Gilain is referring to Cauchy's Collected Works, S\'erie 2,
Tome~3, p.~19, corresponding to (Cauchy 1821, \cite[p.~4]{Ca21}).
Both of Gilain's claims are erroneous, as we show below.  Already in
1973, Hourya Benis Sinaceur warned:
\begin{quote}
on dit trop rapidement que c'est Cauchy qui a introduit la `m\'ethode
des limites' entendant par l\`a, plus ou moins vaguement, l'emploi
syst\'ematique de l'{\em \'epsilonisation\/} (Si\-na\-ceur 1973,
\cite[p.~108]{Si}).
\end{quote}
Sinaceur pointed out that Cauchy's definition of limit resembles, not
that of Weierstrass, but rather that of Lacroix dating from 1810 (see
Sinaceur \cite[p.~109]{Si}).

Cauchy's primary notion is that of {\em variable quantities\/},
introduced in the following terms:
\begin{quote}
On nomme quantit\'e {\em variable\/} celle que l'on consid\`ere comme
devant re\c cevoir successivement plusieurs valeurs diff\'erentes les
unes des autres (Cauchy 1821, \cite[p.~4]{Ca21}).
\end{quote}
Next, Cauchy exploits his primary notion to evoke his kinetic concept
of limit in the following terms:
\begin{quote}
Lorsque les valeurs successivement attribu\'ees \`a une m\^eme
variable s'approchent ind\'efiniment d'une valeur fixe, de mani\`ere
\`a finir par en diff\'erer par aussi peu qu l'on voudra, cette
derni\`ere est appel\'ee la {\em limite\/} de toutes les autres
(ibid.).
\end{quote}
Finally, Cauchy proceeds to define infinitesimals in the following
terms: 
\begin{quote}
Lorsque les valeurs num\'eriques successives d'une m\^eme variable
d\'ecroissent ind\'efiniment, de mani\`ere \`a s'abaisser au-dessous
de tout nombre donn\'e, cette variable devient ce qu'on nomme {\em un
infiniment petit\/} ou une quantit\'e infiniment petite.  Une variable
de cette esp\`ece a z\'ero pour limite (ibid.).
\end{quote}
Thus, Cauchy defined both infinitesimals and limits in terms of
variable quantities.  Neither is the limit concept primary, nor are
infinitesimals defined in terms of limits, contrary to Gilain's
claims.

\subsection{Hawking and Gray}

In his 2007 anthology, S.~Hawking reproduces Cauchy's {\em
infinitesimal\/} definition of continuity on page 639:
\begin{quote}
the function~$f(x)$ remains continuous with respect to~$x$ between the
given bounds,%
\footnote{\label{f1}The term \emph{bounds} is Hawking's translation of
Cauchy's \emph{limites}.}
if, between these bounds, an infinitely small increment in the
variable always produces an infinitely small increment in the function
itself (Hawking, \cite[p.~639]{Ha}).
\end{quote}
But Hawking also claims \emph{on the same page} 639, in a comic {\em
non-sequitur\/}, that Cauchy ``was particularly concerned to banish
infinitesimals''.

In a similar vein, historian J.~Gray lists \emph{continuity} among
concepts Cauchy allegedly defined
\begin{quote}
using careful, if not altogether unambiguous, {\bf limiting} arguments
(Gray 2008, \cite[p.~62]{Gray08}) [emphasis added--authors].
\end{quote}
But in point of fact {\em limits\/} appear in Cauchy's definition only
in the sense of the \emph{extremities} (endpoints, bounds) of the
domain of definition.%
\footnote{See footnote~\ref{f1}.}
Contrary to Gray's claim, the arguments Cauchy used to define
continuity were not ``limiting'' but rather infinitesimal.

\subsection{Laugwitz replies}

Dieudonn\'e's query mentioned in Subsection~\ref{12} is answered by
Laugwitz,%
\footnote{Laugwitz reports having been influenced by Freudenthal's
allusion to Cauchy's work on delta functions in \cite[p.~136]{Fr}.}
who argues that objects such as delta functions (and their potential
applications) disappeared from the literature due to the elimination
of infinitesimals, in whose absence they could not be sustained.
Laugwitz notes that
\begin{quote}
Cauchy's use of delta function methods in Fourier analysis and in the
summation of divergent integrals enables us to analyze the change of
his attitude toward infinitesimals (Laugwitz 1989,
\cite[p.~232]{Lau89}).
\end{quote}

\section{From Cauchy to Dirac}

A function of the type generally attributed to Dirac (1902--1984) was
specifically described in 1827 by Cauchy in terms of infinitesimals.
More specifically, Cauchy uses a unit-impulse, infinitely tall,
infinitely narrow delta function, as an integral kernel.  Thus, in
1827, Cauchy used infinitesimals in his definition of a ``Dirac''
delta function \cite[p.~188]{Ca27}.  Here Cauchy uses
infinitesimals~$\alpha$ and~$\epsilon$, where~$\alpha$ is, in modern
terms, the ``scale parameter'' of the ``Cauchy distribution'', whereas
$\epsilon$ gives the size of the interval of integration.  Cauchy
wrote \cite[p.~188]{Ca27}:
\begin{quote}
Moreover one finds, denoting by~$\alpha$,~$\epsilon$ two infinitely
small numbers,
\begin{equation}
\label{151}
\frac{1}{2} \int_{a-\epsilon}^{a+\epsilon} F(\mu) \frac{\alpha \;
d\mu}{\alpha^2 + (\mu-a)^2} = \frac{\pi}{2} F(a)
\end{equation}
\end{quote}
(Cauchy's 1815-1827 text is analyzed in more detail in
Section~\ref{XVIII}).%
\footnote{See Laugwitz \cite[p.~230]{Lau89}.}
A formula equivalent to~\eqref{151} was proposed by Dirac a century
later.%
\footnote{See main text at footnote~\ref{f3}.}
The expression
\begin{equation}
\label{31b}
\frac{\alpha}{\alpha^2 + (\mu-a)^2}
\end{equation}
(for real~$\alpha$) is known as the {\em Cauchy distribution\/} in
probability theory.  Here Cauchy specifies a function which meets the
criteria as set forth by Dirac a century later.%
\footnote{\label{f7}See footnote~\ref{spec1} for a discussion of
Dirac's criterion.}
In modern terminology, the function is called the probability density
function, and the parameter~$\alpha$ is referred to as the {\em scale
parameter\/}.  Cauchy integrates~$F$ against the kernel~\eqref{31b} as
in~\eqref{151} so as to extract the value of~$F$ at the point~$a$,
exploiting the characteristic property of a delta function.%
\footnote{\label{ff}From the modern viewpoint, formula~\eqref{151}
holds up to an infinitesimal error; thus, the right hand side is the
standard part of the left hand side, see Section~\ref{rival}.}
Thus, a Cauchy distribution with an infinitesimal scale parameter
produces an entity with Dirac-delta function behavior, exploited by
Cauchy already in 1827.

Laugwitz notes that Cauchy's formula \eqref{151} is satisfied when
$\epsilon\geq \alpha^{1/2}$ (as well as for all positive real values
of~$\epsilon>0$).  Today we would write the integral in the formally
more elegant fashion as an integral over an infinite domain:%
\footnote{Modulo the proviso of footnote~\ref{ff}.}
\[
\frac{1}{2} \int_{-\infty}^{\infty} F(\mu) \frac{\alpha \;
d\mu}{\alpha^2 + (\mu-a)^2} = \frac{\pi}{2} F(a),
\]
but Cauchy, ever the practical man rather than a formalist, saw no
reason to bother with subjunctive integration over domains where the
function is perceptually indistinguishable from zero anyway.

Furthermore, Laugwitz documents Cauchy's use of
\begin{quote}
an explicit delta function \emph{not} contained under an integral sign
(Laugwitz 1989, \cite[p.~231]{Lau89}) [emphasis added---authors],
\end{quote}
contrary to a claim in Dieudonn\'e's text.%
\footnote{\label{notunder1}See discussion of Dieudonn\'e's claim in
footnote~\ref{notunder2}.}
An occurrence of a delta function {\em not\/} under an integral sign
in Cauchy's work is discussed in Section~\ref{notunder3}.

In his 1908 book \emph{Elementary Mathematics from an Advanced
Standpoint}, Felix Klein points out that the
\begin{quote}
na\"\i ve [perceptual] methods always rise to unconscious importance
whenever in mathematical physics, mechanics, or differential geometry
a preliminary theorem is to be set up.  You all know that they are
very serviceable then (Klein \cite[p.~211]{Kl08}).
\end{quote}
These remarks apply perfectly well to Cauchy's~$\delta$, as well.  On
the other hand, Klein is perfectly aware of the situation on the
ground:
\begin{quote}
To be sure, the pure mathematician is not sparing of his scorn on
these occasions.  When I was a student it was said that the
differential, for a physicist, was a piece of brass which he treated
as he did the rest of his apparatus (ibid).
\end{quote}
Additional remarks by Klein, showing the importance he attached to
this vital connection, appear in Section~\ref{klein}.

\section{Cauchy's Note XVIII}
\label{XVIII}

Cauchy's lengthy work {\em Th\'eorie de la propagation des ondes \`a
la surface d'un fluide pesant d'une profondeur ind\'efinie}
\cite{Ca27} was written in 1815.  The manuscript was published in 1827
as a 300-page text, with a number of additional Notes at the end.  The
running shortened title used throughout is {\em M\'emoire sur la
th\'eorie des ondes\/}.

Note XVIII, entitled {\em Sur les int\'egrales d\'efinies
singuli\`eres et les valeurs principales des integrales
indetermin\'ees\/}, starts on page 288.  Cauchy recalls the notion of
a singular definite integral, describing it in terms of an integrand
that becomes ``infinite or indeterminate''.  He continues by denoting
by~$\varepsilon$ an ``infinitely small number''%
\footnote{Cauchy's use of term ``number'', rather than ``quantity'',
in this context is interesting.  In general he refers to
infinitesimals as ``quantities''.}
and by~$a,b$ two positive constants.  On page 289, after choosing an
additional ``infinitely small number''~$\alpha$, Cauchy writes down
the integral
\begin{equation*}
\frac{1}{2} \int_{a-\epsilon}^{a+\epsilon} F(\mu) \frac{\alpha \;
d\mu}{\alpha^2 + (\mu-a)^2} = \frac{\pi}{2} F(a)
\end{equation*}
(already reproduced as formula~\eqref{151} above), which he denotes
by~(2).  Cauchy proceeds to point out that, since the integrand of his
equation~(2) is {\em sensiblement \'egale \`a z\'ero\/} [essentially
equal to zero] for all values of~$\mu$ {\em qui ne sont pas tr\`es
rapproch\'ees de~$a$} [which are not too close to~$a$], it follows
that the integrals appearing in his earlier Note VI reduce to singular
integrals determined by his equation (2).

Note XVIII then proceeds to discuss principal values and to offer
alternative derivations of a number of earlier results, and is
concluded on page 299.

\section{Cauchy's 1827 {\em M\'emoire\/}}
\label{notunder3}

An additional occurrence of a delta function occurs in Cauchy's brief
1827 text {\em M\'emoire sur les d\'eveloppements des fonctions en
s\'eries p\'eriodiques\/}~\cite{Ca27a}.  The text contains an ({\em a
priori\/} doomed) attempt to prove the convergence of Fourier series
under the sole assumption of continuity.  What concerns us here is
his, correct, use of infinitesimals at a certain stage in the
argument.  Cauchy opens his {\em m\'emoire\/} with a discussion of the
importance of what are known today as Fourier series, in ``a large
number of problems of mathematical physics'' \cite[p.~12]{Ca27a}.  On
page 13, Cauchy denotes by~$\varepsilon$ {\em un nombre infiniment
petit\/} [an infinitely small number], lets~$\theta= 1-\varepsilon$,
and lets~$x$ be between~$0$ and~$a=2\pi$.  On page~14, he points out
that the expression
\begin{equation}
\label{14}
1 + \frac{1}{e^{-i(x-\mu)}-\theta} + \frac{1}{e^{i(x-\mu)}-\theta}
\end{equation}
(his notation is slightly different) ``will be essentially zero,
except when~$\mu$ differs very little from~$x$''.  Note that the
expression \eqref{14} appearing on Cauchy's page 14, does {\em not\/}
occur under the integral sign (it was exploited as a kernel in the
last formula on the previous page~13).  

Cauchy then sets~$\mu =x+iw$ and concludes that the integral will be
essentially reduced to%
\footnote{Cauchy wrote~$a$ where we wrote~$2\pi$.}
\begin{equation*}
f(x) \cdot \int_{\frac{-x}{\varepsilon}}^{\frac{2\pi-x}{\varepsilon}}
\left( \frac{1}{1+iw} + \frac{1}{1-iw} \right) dw = 2\pi f(x).
\end{equation*}

\section{Triumvirate history}
\label{four}

Studies seeking to document continuity between Cauchy's infinitesimals
and modern set-theoretic implementations of infinitesimals have not
always been viewed sympathetically by commentators.

In 1949, C.~Boyer hagiographically described G.~Cantor, R.~Dedekind,
and K.~Weierstrass as ``the great triumvirate'' \cite[p.~298]{Boy}.
The ``triumvirate'' historian tends to view the history of mathematics
as an ineluctable march toward the radiant future of Weierstrassian
epsilontics.  Such a stance leaves little room for a sympathetic view
of a continuity between infinitesimals as practiced prior to the
triumvirate, on the one hand, and 20th century set-theoretic
\emph{implementations} of infinitesimals, on the other.  J.~Dauben
wrote:
\begin{quote}
In bringing historians of the calculus to task, Robinson was
particularly critical of Carl Boyer's \emph{The Concepts of the
Calculus}.  The history of mathematics was on shaky ground, Robinson
felt, if it chose to pass judgment on earlier theories based upon
currently fashionable prejudices.  Nonstandard analysis cast a new
light on the history of the calculus, and Robinson was interested to
see how it might appear if reexamined without assuming that
infinitesimals were wrongheaded or at all lacking in rigor (Dauben
1995, \cite[p.~349]{Da95}).
\end{quote}
Dauben appears sympathetic to Robinson's interest in re-examining the
history of infinitesimals.  It is all the more surprising therefore to
read a text Dauben authored a few years earlier.  Here Dauben presents
a list of authors, including D.~Laugwitz, who ``have used nonstandard
analysis to rehabilitate or `vindicate' earlier infinitesimalists'',
and concludes:
\begin{quote}
Leibniz, Euler, and Cauchy [\ldots] had, in the views of some
commentators, ``Robinsonian'' nonstandard infinitesimals in mind from
the beginning.  The most detailed and methodically [sic] sophisticated
of such treatments to date is that provided by Imre Lakatos; in what
follows, it is his analysis of Cauchy that is emphasized (Dauben 1988,
\cite[p.~179]{Da88}).
\end{quote}
However, Lakatos's treatment was certainly not ``the most detailed and
methodically sophisticated'' one by the time Dauben's text appeared in
1988.  Thus, in 1987, Laugwitz had published a detailed scholarly
study of Cauchy in \emph{Historia Mathematica} (Laugwitz
\cite{Lau87}).  Laugwitz's text in \emph{Historia Mathematica} appears
to be the published version of his 1985 preprint \emph{Cauchy and
infinitesimals}.  Laugwitz's 1985 preprint does appear in Dauben's
bibliography (Dauben, 1988 \cite[p.~199]{Da88}), indicating that
Dauben was familiar with it.  It is odd to suggest, as Dauben seems
to, that a scholarly study published in \emph{Historia Mathematica}
would countenance a view that Cauchy could have had ```Robinsonian'
nonstandard infinitesimals in mind from the beginning''.  Surely
Dauben has committed a strawman fallacy here.

Rather, Lakatos, Laugwitz, Cutland et al. \cite{CKKR}, Br\aa ting
(2007, \cite{Br}) and others have argued that infinitesimals as
employed by Cauchy have found set-theoretic implementations in the
framework of modern theories of infinitesimals, just as Kanovei had
done for Euler in 1988 \cite{Kan}.  The existence of such
implementations indicates that the historical infinitesimals were less
prone to contradiction than has been routinely maintained by
triumvirate historians, who invariably cite Berkeley's flawed
empiricist critique.  The issue is dealt with in more detail by Katz
\& Katz \cite{KK11a}, \cite{KK11b}, \cite{KK11c}, \cite{KK11d}; B\l
aszczyk et al. \cite{BKS}; Borovik et al.~\cite{BK}; Katz \&
Leichtnam~\cite{KL}; and Katz \& Sherry \cite{KS1, KS2}.

What criteria can we employ to evaluate the achievements of modern
theories of infinitesimals?  Remarkably, just such a criterion was
explicitly provided by both Felix Klein and Abraham Fraenkel, as
discussed in Section~\ref{KF}.

\section{Klein--Fraenkel criterion}
\label{KF}

In 1908, Klein formulated a criterion of what it would take for a
theory of infinitesimals to be successful.  Namely, one must be able
to prove a mean value theorem for arbitrary intervals, including
infinitesimal ones:
\begin{quote}
The question naturally arises whether [\ldots] it would be possible to
modify the traditional foundations of infinitesimal calculus, so as to
include actually infinitely small quantities in a way that would
satisfy modern demands as to rigor; in other words, to construct a
non-Archimedean system.  The first and chief problem of this analysis
would be to prove the mean-value theorem
\[
f(x+h)-f(x)=h \cdot f'(x+\vartheta h)
\]
from the assumed axioms.  I will not say that progress in this
direction is impossible, but it is true that none of the investigators
have achieved anything positive (Klein 1908, \cite[p.~219]{Kl08}).
\end{quote}
In 1928, A.~Fraenkel \cite[pp.~116-117]{Fran} formulated a similar
requirement in terms of the mean value theorem.  Such a
Klein--Fraenkel criterion is satisfied by the Hewitt-\Los-Robinson
theory.  Indeed, the mean value theorem is true for an arbitrary
hyperreal interval by the transfer principle.  Fraenkel's opinion of
Robinson's theory is on record:
\begin{quote}
my former student Abraham Robinson had succeeded in saving the honour
of infinitesimals - although in quite a different way than Cohen%
\footnote{The reference is to Hermann Cohen (1842--1918), whose
fascination with infinitesimals elicited fierce criticism by both
Cantor and B.~Russell.  For an analysis of Russell's
\emph{non-sequiturs}, see Ehrlich \cite{Eh06} and Katz \& Sherry
\cite{KS1}.}
and his school had imagined (Fraenkel 1967, \cite[p.~107]{Fra67}).
\end{quote}

\section{Set-theoretic implementation of infinitesimals}
\label{rival}

In Section~\ref{four}, we clarified the role of modern theories of
infinitesimals in interpreting the work of historical
infinitesimalists.  Here we present some details of a particular
set-theoretic implementation of infinitesimals as developed through
the work of Hewitt, \Los, and Robinson.  For an alternative approach
to infinitesimals, see P.~Giordano~\cite{Gio10b}.

In 1961, Robinson \cite{Ro61} constructed an infinitesimal-enriched
continuum, suitable for use in calculus, analysis, and elsewhere,
based on earlier work by E.~Hewitt \cite{Hew}, J.~\Los{} \cite{Lo},
and others.  In 1962, W.~Luxemburg \cite{Lu62} popularized a
presentation of Robinson's theory in terms of the ultrapower
construction,%
\footnote{Note that both the term ``hyper-real'', and an ultrapower
construction of a hyperreal field, are due to E.~Hewitt in 1948, see
\cite[p.~74]{Hew}.  Luxemburg~\cite{Lu62} clarified its relation to
the competing construction of Schmieden and Laugwitz \cite{SL}, also
based on sequences, but using the ideal~$\fez$.  Dauben
\cite[p.~395]{Da95} mistakenly suggests that it was Luxemburg who
initiated the ultrapower approach to the hyperreals using free
ultrafilters.}
in the mainstream foundational framework of the Zermelo-Fraenkel set
theory with the axiom of choice.  Namely, the hyperreal field is the
quotient of the collection of arbitrary sequences, where a sequence
\begin{equation}
\label{11}
\langle u_1, u_2, u_3, \ldots \rangle
\end{equation}
converging to zero generates an infinitesimal.  Arithmetic operations
are defined at the level of representing sequences; e.g., addition and
multiplication are defined term-by-term.

To motivate the construction of the hyperreals, we will start with the
ring~$\Q^\N$ of sequences of rational numbers.  Let~$\mathcal{C}_\Q
\subset \Q^\N$ denote the subring consisting of Cauchy sequences.  The
reals are by definition the quotient field%
\footnote{\label{meray}Such a construction is usually attributed to
Cantor, and is actually due to M\'eray 1869, \cite{Me} who published
three years earlier than E.~Heine.}
\begin{equation}
\label{realbis}
\R:= \mathcal{C}_\Q / \fnull,
\end{equation}
where~$\fnull$ contains all null sequences.  An infinitesimal-enriched
extension of~$\Q$ may be obtained by modifying~\eqref{realbis} as
follows.  We consider a subring~$\fez\subset\fnull$ of sequences that
are ``eventually zero'', i.e., vanish at all but finitely many places.
Then the quotient~$\mathcal{C}_\Q / \fez$ naturally surjects onto~$\R=
\mathcal{C}_\Q / \fnull$.  The elements in the kernel of the
surjection~$\mathcal{C}_\Q/\fez \to \R$ are prototypes of
infinitesimals.  Note that the quotient~$\mathcal{C}_\Q/\fez$ is not a
field, as~$\fez$ is not a maximal ideal.  It is more convenient to
describe the modified construction using the ring~$\R^\N$ rather than
$\mathcal{C}_\Q$.

We therefore redefine~$\fez$ to be the ring of sequences in~$\R^\N$
that eventually vanish, and choose a \emph{maximal} proper ideal
$\mathcal{M}$ so that we have
\begin{equation}
\label{23}
\fez\subset\mathcal{M}\subset\R^\N.
\end{equation}
Then the quotient~$\RRR:=\R^\N/\mathcal{M}$ is a hyperreal field.  The
foundational material needed to ensure the existence of a maximal
ideal~$\mathcal{M}$ satisfing~\eqref{23} is weaker than the axiom of
choice.  This concludes the construction of a hyperreal field~$\RRR$
in the traditional foundational framework, Zermelo-Fraenkel set theory
with the axiom of choice (ZFC).

Let~$I\subset\RRR$ be the subring consisting of infinitesimal elements
(i.e., elements~$e$ such that~$|e|<\frac{1}{n}$ for all~$n\in\N$).
Denote by~$I^{-1}$ the set of inverses of nonzero elements of~$I$.
The complement~$\RRR\setminus I^{-1}$ consists of all the finite
(sometimes called \emph{limited}) hyperreals.  Constant sequences
provide an inclusion~$\R\subset\RRR$.  Every element~$x\in
\RRR\setminus I^{-1}$ is infinitely close to some real number
$x_0\in\R$.  The \emph{standard part function}, denoted ``st'',
associates to every finite hyperreal, the unique real infinitely close
to it:
\[
\text{st}:\RRR\setminus I^{-1} \to \R, \quad x \mapsto x_0.
\]
If~$x$ happens to be represented by a Cauchy sequence~$\langle x_n :
n\in\N \rangle$, so that~$x=\left[ \langle x_n \rangle \right]$, then
the standard part can be expressed in terms of the ordinary limit:
\[
\text{st}(x)=\lim_{n\to \infty} x_n.
\]
More advanced properties of the hyperreals such as saturation were
proved later (see Keisler \cite{Ke94} for a historical outline).  A
helpful ``semicolon'' notation for presenting an extended decimal
expansion of a hyperreal was described by A.~H.~Lightstone~\cite{Li}.
See also P.~Roquette~\cite{Roq} for infinitesimal reminiscences.  A
discussion of infinitesimal optics is in K.~Stroyan \cite{Str},
J.~Keisler~\cite{Ke} and others.  Applications of
infinitesimal-enriched continua range from aid in teaching calculus
\cite{El, KK1, KK2} to the Bolzmann equation (see
L.~Arkeryd~\cite{Ar81, Ar05}) and mathematical physics (see Albeverio
et al. \cite{Alb}).  Edward Nelson \cite{Ne} in 1977 proposed an
alternative to ZFC which is a richer (more stratified) axiomatisation
for set theory, called Internal Set Theory (IST), more congenial to
infinitesimals than ZFC.  The traditional construction of the reals
out of Cauchy sequences can be factored through the hyperreals (see
Giordano et al. \cite{GK11}).  The hyperreals can be constructed out
of integers (see Borovik et al.~\cite{BJK}).

\section{Heaviside function}
\label{nine}

Having clarified the connection between historical infinitesimals and
their modern set-theoretic implementations in Sections~\ref{four}
and~\ref{rival}, we return to delta functions and their integrals.

Dieudonn\'e's review of L\"utzen's book is assorted with the habitual,
and near-ritual on the part of some mathematicians, expression of
disdain for physicists:
\begin{quote}
However, a function such as the Heaviside function on~$\R$, equal to
$1$ for~$x\geq 0$ and to~$0$ for~$x<0$, has no weak derivative, in
spite of its very mild discontinuity; at least this is what the
mathematicians would say, but physicists thought otherwise, since for
them there {\em was\/} a ``derivative''~$\delta$, the Dirac ``delta
function'' (Dieudonn\'e \cite[p.~377]{Die}) [the quotation marks are
Dieudonn\'e's---authors].
\end{quote}
Dieudonn\'e then proceeds to make the following remarkable claims:
\begin{quote}
Of course, there was before 1936 no reasonable mathematical definition
of these objects; but it is characteristic that they were never used
in {\em bona fide\/} computations except {\em under the integral
sign\/},%
\footnote{\label{notunder2}We are not certain what {\em bona fide\/}
calculations are exactly, but at any rate Dieudonn\'e's claim that
delta functions were never used except under an integral sign, would
be inaccurate; see main text at footnote~\ref{notunder1}.}
giving formulas%
\footnote{Here we have simplified Dieudonn\'e's formula in
\cite[p.~377]{Die}, by restricting to the special case~$n=0$.}
such as
\begin{equation*}
\int \delta(x-a) f(x) = f(a).
\end{equation*}
\end{quote}
Are Dieudonn\'e's claims accurate?  Dieudonn\'e's claim that, before
1936, delta functions occurred only under the integral sign, is
contradicted by Cauchy's use of a delta function {\em not\/} contained
under an integral sign, over a hundred years earlier (see
Section~\ref{notunder3}).

Are the physicists so far off the mark in speaking of the delta
function as the derivative of the Heaviside function?  Is it really
true that there was no reasonable mathematical definition before 1936,
as Dieudonn\'e claims?  In fact, Cauchy had a reasonable mathematical
definition of a delta function, though of course both set theory and
set-theoretic implementations of his ideas were still decades away.

\begin{figure}
\includegraphics[height=2in]{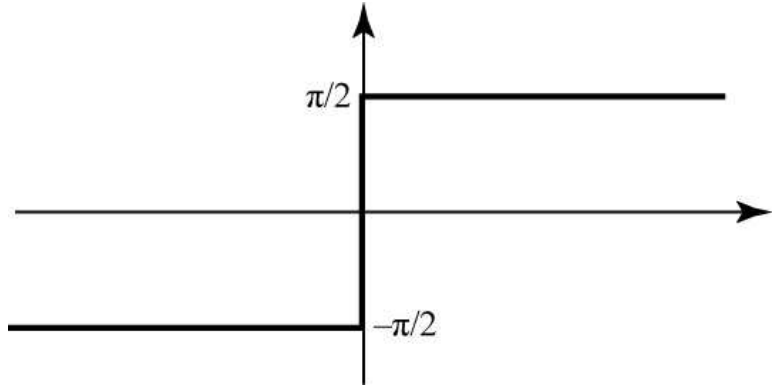}
\caption{Heaviside function}
\label{heavy}
\end{figure}

Thus, consider the zigzag~$\mathcal{Z}\subset \R^2$ in
the~$(x,y)$-plane given by the union
\begin{equation}
\label{zig}
\mathcal{Z}= \left( \R_- \times \left\{ -\tfrac{\pi}{2} \right\}
\right) \cup \left( \{0\} \times \left[ -\tfrac{\pi}{2} \tfrac{\pi}{2}
\right] \right) \cup \left( \R_+ \times \left\{ + \tfrac{\pi}{2}
\right\} \right),
\end{equation}
thought of as the physicist's Heaviside function, see
Figure~\ref{heavy}.

Now consider the graph of~$\arctan(x)$ in the~$(x,y)$-plane, and
compress it toward the~$y$-axis by means of a sequence of
functions~$\arctan(nx)$, or~$\arctan(x/\alpha)$
where~$\alpha=\frac{1}{n}$.  Their derivatives~$F_\alpha(x)$ satisfy
\begin{equation*}
\int_{-\infty}^\infty F_\alpha = \pi
\end{equation*}
by the fundamental theorem of calculus.  In an infinitesimal-enriched
continuum (see Section~\ref{rival}), we can assign an infinitesimal
value to~$\alpha$.  Then the graph
\begin{equation*}
\Gamma_{\arctan(x/\alpha)}
\end{equation*}
of~$\arctan(x/\alpha)$ is ``appreciably indistinguishable'' from the
plane zigzag
$\mathcal{Z}\subset\R^2$ of formula~\eqref{zig}.%
\footnote{In modern notation, this relation would be expressed by the
fact that the standard part ``st'' of the graph
$\Gamma_{\arctan(x/\alpha)}$ coincides with the zigzag:
\[
{\rm st} \left( \Gamma_{\arctan(x/\alpha)} \right) = \mathcal{Z}
\subset \R^2.
\]
Here the function~$\arctan(x/\alpha)$ is the mathematical counterpart
of the physicist's Heaviside function.  Of course, Cauchy did not have
the notion of a standard part function (see Section~\ref{rival}), to
express the idea that an error term is infinitesimal.  Instead, he
used the expression {\em sensiblement nulle\/} (sensibly nothing), see
\cite[p.~231]{Lau89}.}
Instead of attempting to differentiate the zigzag itself as physicists
are alleged to do, we differentiate its infinitesimal
approximation~$\arctan(x/\alpha)$, and note that we obtain precisely
Cauchy's delta function appearing in formula~\eqref{151}, against
which~$F$ is integrated.%
\footnote{Takeuti \cite{Ta}, Giorello \cite{Gio}, Lightstone and Wong
\cite{LW}, and later Todorov~\cite{To} and P\'eraire~\cite{Per} have
developed this theme further.  Yamashita~\cite{Ya} provides a
bibliography of articles dealing with hyperreal delta functions.}

\section{Klein's remarks on physics}
\label{klein}

As in the case of Cauchy's delta function, infinitesimals provide an
intuitive point of entry to key phenomena in both mathematics and
physics.  In a similar vein, Klein discussed infinitesimal
oscillations of the pendulum in his \emph{Elementary Mathematics from
an Advanced Standpoint}.  Klein presents the derivation of the
pendulum law by pointing out that
\begin{quote}
it follows from the fundamental laws of mechanics that the motion of
the pendulum is determined by the equation~$\frac{d^2\phi}{dt^2} = \
\frac{g}{\ell} \sin \phi$ (Klein \cite[p.~187]{Kl08}).
\end{quote}
Here~$g$ is the gravitational constant, while~$\ell$ is the length of
the thread by which the pendulum is suspended, and~$\phi$ is the angle
of deviation from the normal.  Klein continues:
\begin{quote}
For small amplitudes we may replace~$\sin \phi$ by~$\phi$ without
serious error.  This gives for the so called infinitely small
oscillation of the pendulum~$\frac{d^2\phi}{dt^2} = \ \frac{g}{\ell}
\phi$ (ibid.).
\end{quote}
Klein proceeds to write down the general solution~$\phi = C \cos
\sqrt{\tfrac{g}{\ell}} (t-t_0)$, and points out that the duration of a
complete oscillation, i.e., the period~$T=2\pi\sqrt{\ell/g}$, is
independent of the amplitude~$C$.  Reflecting upon the teaching
practices at the time, Klein muses over the incongruity of
\begin{quote}
the curious phenomenon that one and the same teacher, during one hour,
the one devoted to mathematics, makes the very highest demands as to
the logical exactness of all conclusions.  In his judgment [\ldots]
his demands are not satisfied by the infinitesimal calculus.  In the
next hour, however, that devoted to physics, he accepts the most
questionable conclusions and makes the most daring applications of
infinitesimals (ibid.).
\end{quote}
Klein's lament concerning the existence of such artificial boundaries
is echoed by Novikov \cite{No2}.

\section{Charles Sanders Peirce's framework}

The customary set-theoretic framework that has become the reflexive
litmus test of {\em mathematical rigor\/} in most fields of modern
mathematics (with the possible exception of the field of mathematical
logic) makes it difficult to analyze Cauchy's use of infinitesimals,
and to evaluate its significance.  We will therefore use a conceptual
framework proposed by C. S. Peirce in 1897, in the context of his
analysis of the concept of continuity and continuum, which, as he felt
at the time, is composed of infinitesimal parts, see
\cite[p.~103]{Ha}.  Peirce identified three stages in creating a novel
concept:

\begin{quote}
there are three grades of clearness in our apprehensions of the
meanings of words.  The first consists in the connexion of the word
with familiar experience. \ldots The second grade consists in the
abstract definition, depending upon an analysis of just what it is
that makes the word applicable. \ldots The third grade of clearness
consists in such a representation of the idea that fruitful reasoning
can be made to turn upon it, and that it can be applied to the
resolution of difficult practical problems \cite{Pei} (see Havenel
2008, \cite[p.~87]{Ha}).
\end{quote}

The ``three grades'' can therefore be summarized as
\begin{enumerate}
\item
familiarity through experience;
\item
abstract definition aimed at applications;
\item
fruitful reasoning ``made to turn'' upon it, with applications.
\end{enumerate}

To apply Peirce's framework to Cauchy's notion of infinitesimal, we
note that grade (1) is captured in Cauchy's description of continuity
of a function in terms of ``varying by imperceptible degrees''.  Such
a turn of phrase occurs both in his letter to Coriolis (see Cauchy
1837, \cite{Ca37}), and in his 1853 text \cite[p.~35]{Ca53}.%
\footnote{Note that both Cauchy's original French ``par degr\'es
insensibles'', and its correct English translation ``by imperceptible
degrees'', are etymologically related to {\em sensory perception\/}.}
At Grade (2), Cauchy describes infinitesimals as generated by null
sequences (see \cite{Br}), and defines continuity in terms of an
infinitesimal~$x$-increment resulting in an infinitesimal change
in~$y$.  Finally, at stage (3), Cauchy fruitfully applies the
crystallized notion of infinitesimal both in Fourier analysis and in
evaluation of singular integrals, by means of a ``Dirac'' delta
function defined in terms of a (Cauchy) distribution with an
infinitesimal scaling parameter.

From the viewpoint of the Peircian framework, and contrary to
Dieudonn\'e's claim, Cauchy did have a \emph{reasonable mathematical
definition} of a, Dirac, delta function, though both set-theory and
set-theoretic implementations of infinitesimal-enriched continua in
which Cauchy's definition could be made operative were still decades
away.  Cauchy's definition could be compared to the definitions of the
fundamental concepts of infinitesimal calculus furnished by Newton and
Leibniz.  The founders of the calculus, similarly to Cauchy, lacked a
set-theoretic formalisation of a continuum, and yet are (rightfully)
given credit for the fundamental concepts they introduced.  Cauchy
exploited infinitesimals both in his definition of continuity in 1821,
and in his definition of a notion closely related to uniform
convergence in 1853 (see Katz \& Katz \cite{KK11b} and B\l aszczyk et
al.~\cite{BKS}).  The centrality of infinitesimals in Cauchy's
approach to analysis is further clarified through his use thereof in
defining ``Dirac'' delta functions.

\section*{Acknowledgments}

We are grateful to S. Kutateladze and D. Ross for insightful comments
that helped improve an earlier version of the manuscript.  The
influence of Hilton Kramer (1928-2012) is obvious.

\bigskip

\textbf{Mikhail G. Katz} is Professor of Mathematics at Bar Ilan
University, Ramat Gan, Israel.  Among his publications are the
following: (with P.~B\l aszczyk and D.~Sherry) Ten misconceptions from
the history of analysis and their debunking, \emph{Foundations of
Science} (online first); (with A.~Borovik) Who gave you the
Cauchy--Weierstras tale?  The dual history of rigorous calculus,
\emph{Foundations of Science} \textbf{17} (2012), no.~3, 245-276;
(with A.~Borovik and R.~Jin) An integer construction of
infinitesimals: Toward a theory of Eudoxus hyperreals, \emph{Notre
Dame Journal of Formal Logic} \textbf{53} (2012), no.~4, 557-570;
(with K.~Katz) Cauchy's continuum, \emph{Perspectives on Science}
\textbf{19} (2011), no.~4, 426-452; (with K.~Katz) Meaning in
classical mathematics: is it at odds with Intuitionism?
\emph{Intellectica} \textbf{56} (2011), no.~2, 223-302; (with K.~Katz)
A~Burgessian critique of nominalistic tendencies in contemporary
mathematics and its historiography, \emph{Foundations of Science}
\textbf{17} (2012), no.~1, 51--89; (with K.~Katz) Stevin numbers and
reality, \emph{Foundations of Science} \textbf{17} (2012), no.~2,
109-123; (with E.~Leichtnam) Commuting and non-commuting
infinitesimals, to appear in \emph{American Mathematical Monthly};
(with D.~Schaps and S.~Shnider) Almost Equal: The Method of Adequality
from Diophantus to Fermat and Beyond, \emph{Perspectives on Science}
\textbf{20} (2012), to appear; (with D.~Sherry) Leibniz's
infinitesimals: Their fictionality, their modern implementations, and
their foes from Berkeley to Russell and beyond, \emph{Erkenntnis}
(online first); (with D.~Sherry) Leibniz's laws of continuity and
homogeneity, \emph{Notices of the American Mathematical Society}
\textbf{59} (2012), to appear.

\bigskip

\textbf{David Tall} is Professor Emeritus in the Department of
Education, University of Warwick, CV4 7AL.  He studied for a doctorate
in mathematics with Fields Medallist Michael Atiyah and a doctorate in
Education with Richard Skemp, author of \emph{The Psychology of
Learning Mathematics}.  He has made a life-long study of the
development of mathematical thinking from child to adult, introducing
new concepts, such as the distinction between concept definition and
concept definition with Shlomo Vinner which was selected by the
National Council of Teachers in Mathematics as one of the seventeen
papers of the twentieth century that should be widely read. He
programmed \emph{A Graphic Approach to the Calculus} in the early
1980s, based on the dynamic idea of 'local straightness'. He edited
the book on \emph{Advanced Mathematical Thinking} in 1990 and
introduced the term `procept' with Eddie Gray in 1991 to represent the
duality, ambiguity and flexibility of symbolism that represents both
process and concept.  Procept and concept image were declared two of
the five major contributions to mathematics education by founder
president, Efraim Fischbein, reviewing the first 25 years of the
International Group for the Psychology of Mathematics Education. He
was the most quoted author in the review of the group's first thirty
years.  In 2004, introduced the theoretical framework of \emph{Three
Worlds of Mathematics}, based on human embodiment (perception and
action), operational symbolism and axiomatic formalism together with
an analysis of the supportive aspects that encourage generalization
and problematic aspects that impede progress. He has cooperated with
Mikhail Katz, applying the three-world framework to the work of Cauchy
and subsequent developments in formalism and the use of
infinitesimals.  His papers are all available from
http://www.warwick.ac.uk/staff/David.Tall.


\begin{thebibliography}{AB}

\bibitem{Alb} Albeverio, S.; H\o egh-Krohn, R.; Fenstad, J.; Lindstr\o
m, T.: Nonstandard methods in stochastic analysis and mathematical
physics.  {\em Pure and Applied Mathematics\/}, \textbf{122}. Academic
Press, Inc., Orlando, FL, 1986.


\bibitem{Ar81} Arkeryd, L.: Intermolecular forces of infinite range
and the Boltzmann equation.  {\em Arch. Rational Mech. Anal.\/}
\textbf{77} (1981), no.~1, 11--21.

\bibitem{Ar05} Arkeryd, L.: Nonstandard analysis.  \emph{American
Mathematical Monthly} \textbf{112} (2005), no.~10, 926-928.


\bibitem{Ba11} Barany, M. J.: revisiting the introduction to Cauchy's
{\em Cours d'analyse\/}.  {\em Historia Mathematica\/} \textbf{38}
(2011), no.~3, 368--388.



\bibitem{BKS} B\l aszczyk, P.; Katz, M.; Sherry, D.: Ten
misconceptions from the history of analysis and their debunking.
\emph{Foundations of Science}, 2012.  See

http://dx.doi.org/10.1007/s10699-012-9285-8

and http://arxiv.org/abs/1202.4153

%
%




\bibitem{BJK} Borovik, A.; Jin, R.; Katz, M.: An integer construction
of infinitesimals: Toward a theory of Eudoxus hyperreals.  \emph{Notre
Dame Journal of Formal Logic} \textbf{53} (2012), no.~4, 557-570.

\bibitem{BK} Borovik, A.; Katz, M.: Who gave you the
Cauchy--Weierstrass tale?  The dual history of rigorous calculus.
\emph{Foundations of Science} \textbf{17} (2012), no.~3, 245-276.  See
http://dx.doi.org/10.1007/s10699-011-9235-x

and http://arxiv.org/abs/1108.2885



\bibitem{Boy} Boyer, C.: The concepts of the calculus, A Critical and
Historical Discussion of the Derivative and the Integral.  Hafner
Publishing Company, 1949.



\bibitem{Br} Br\aa ting, K.: A new look at E. G. Bj\"orling and the
Cauchy sum theorem.  \emph{Archive for History of Exact Sciences}
\textbf{61} (2007), no.~5, 519--535.


\bibitem{Ca21} Cauchy, A. L.: {\em Cours d'Analyse de L'Ecole Royale
Polytechnique.  Premi\`ere Partie.  Analyse alg\'ebrique\/}.  Paris:
Imprim\'erie Royale, 1821.  Online at

{\tiny http://books.google.com/books?id=\_mYVAAAAQAAJ\&dq=cauchy\&lr=\&source=gbs\_navlinks\_s}


\bibitem{Ca27} Cauchy, A.-L. (1815) Th\'eorie de la propagation des
ondes \`a la surface d'un fluide pesant d'une profondeur ind\'efinie
(published 1827, with additional Notes). Oeuvres, Series 1, Vol.~1,
4-318.



\bibitem{Ca27a} Cauchy, A.-L. (1827) M\'emoire sur les
d\'eveloppements des fonctions en s\'eries p\'eriodiques.  Oeuvres,
Series 1, Vol.~2, 12-19.

\bibitem{Ca37} Cauchy, A.-L. (1837) ``Extrait d'une lettre \`a
M. Coriolis," \emph{Oeuvres Compl\`etes}, Series 1, Vol.~4 (Paris:
Gauthier Villars, 1884), 38--42.

\bibitem{Ca53} Cauchy, A. L. (1853) Note sur les s\'eries convergentes
dont les divers termes sont des fonctions continues d'une variable
r\'eelle ou imaginaire, entre des limites donn\'ees. In \emph{Oeuvres
compl\`etes}, Series 1, Vol.~12, pp.~30-36.  Paris: Gauthier--Villars,
1900.


\bibitem{CKKR} Cutland, N.; Kessler, C.; Kopp, E.; Ross, D.: On
Cauchy's notion of infinitesimal.  \emph{The British Journal for the
Philosophy of Science} \textbf{39} (1988), no.~3, 375--378.



\bibitem{Da88} Dauben, J.: Abraham Robinson and Nonstandard Analysis:
History, Philosophy, and Foundations of Mathematics.  In William
Aspray and Philip Kitcher, eds.  History and philosophy of modern
mathematics (Minneapolis, MN, 1985), 177--200, Minnesota
Stud. Philos. Sci., XI, Univ. Minnesota Press, Minneapolis, MN, 1988.
Available online at

http://www.mcps.umn.edu/philosophy/11\_7Dauben.pdf


\bibitem{Da95} Dauben, J.: Abraham Robinson. The creation of
nonstandard analysis.  A personal and mathematical odyssey.  With a
foreword by Benoit B. Mandelbrot.  Princeton University Press,
Princeton, NJ, 1995.





\bibitem{Die} Dieudonn\'e, J.: Reviews: The Prehistory of the Theory of
Distributions.  \emph{American Mathematical Monthly} \textbf{91}
(1984), no.~6, 374--379.

\bibitem{Dir} Dirac, P.: The Principles of Quantum Mechanics.  4th
edition.  Oxford, at the Clarendon Press, 1958.

\bibitem{Eh06} Ehrlich, P.: The rise of non-Archimedean mathematics
and the roots of a misconception. I. The emergence of non-Archimedean
systems of magnitudes.  \emph{Archive for History of Exact Sciences}
\textbf{60} (2006), no.~1, 1--121.


\bibitem{El} Ely, R.: Nonstandard student conceptions about
infinitesimals.  \emph{Journal for Research in Mathematics Education}
\textbf{41} (2010), no.~2, 117-146.


\bibitem{Fi} Fisher, G.: Cauchy's Variables and Orders of the
Infinitely Small.  \emph{The British Journal for the Philosophy of
Science}, \textbf{30} (1979), no.~3, 261--265.

\bibitem{Fi81} Fisher, G.: The infinite and infinitesimal quantities
of du Bois-Reymond and their reception.  \emph{Archive for History of
Exact Sciences} \textbf{24} (1981), no.~2, 101--163.


\bibitem{Fo} Fourier, J.: Th\'eorie analytique de la chaleur.  Paris:
Firmin Didot P\`ere et Fils, 1822.

\bibitem{Fran} Fraenkel, A.: Einleitung in die Mengenlehre.  Dover
Publications, New York, N. Y., 1946 [originally published by Springer,
Berlin, 1928].

\bibitem{Fra67} Fraenkel, A.: Lebenskreise.  Aus den Erinnerungen
eines j\"udischen Mathematikers.  Deutsche Verlags-Anstalt, Stuttgart,
1967.


\bibitem{Fr} Freudenthal, H.: Cauchy, Augustin-Louis.  In Dictionary
of Scientific Biography, ed. by C. C. Gillispie, vol. 3 (New York:
Charles Scribner's sons, 1971), 131-148.  See

http://www.encyclopedia.com/topic/Augustin-Louis\_Cauchy.aspx\#1

\bibitem{Gil} Gilain, C.: Cauchy et le cours d'analyse de l'Ecole
polytechnique.  With an editorial preface by Emmanuel Grison.
\emph{Soci\'et\'e des Amis de la Biblioth\`eque de l'Ecole
Polytechnique.  Bulletin} 1989, no. 5, 145 pp.  

See http://www.sabix.org/bulletin/b5/cauchy.html

\bibitem{Gio10b} Giordano, P.: Infinitesimals without logic.
\emph{Russian Journal of Mathematical Physics}, \emph{17} (2010),
no.~2, 159--191.


\bibitem{GK11} Giordano, P.; Katz, M.: Two ways of obtaining
infinitesimals by refining Cantor's completion of the reals.
Preprint, 2011,

see http://arxiv.org/abs/1109.3553


\bibitem{Gio} Giorello, G.: Una rappresentazione nonstandard delle
distribuzioni temperate e la trasformazione di Fourier.  \emph{Unione
Matematica Italiana.  Bollettino} (4) 7 (1973), 156--167.


\bibitem{Gray08} Gray, J.: Plato's ghost.  The modernist
transformation of mathematics.  Princeton University Press, Princeton,
NJ, 2008.


\bibitem{Ha} Havenel, J.: Peirce's Clarifications of Continuity.
\emph{Transactions of the Charles S. Peirce Society} \textbf{44}
(2008), no.~1, 86-133.

\bibitem{Haw} Hawking, S., Ed.: The mathematical breakthroughs that
changed history.  Running Press, Philadelphia, PA, 2007 (originally
published 2005).

\bibitem{Hew} Hewitt, E.: Rings of real-valued continuous
functions. I.  \emph{Transactions of the American Mathematical
Society} \textbf{64} (1948), 45--99.


\bibitem{Kan} Kanovei, V.: Correctness of the Euler method of
decomposing the sine function into an infinite product. (Russian)
\emph{Uspekhi Mat. Nauk} \textbf{43} (1988), no.~4 (262), 57--81, 255;
translation in \emph{Russian Math. Surveys} \textbf{43} (1988), no.~4,
65--94.


\bibitem{KK1} Katz, K.; Katz, M.: Zooming in on infinitesimal~$1-.9..$
in a post-triumvirate era.  \emph{Educational Studies in Mathematics\/}
\textbf{74} (2010), no.~3, 259-273.  See
{http://arxiv.org/abs/arXiv:1003.1501}



\bibitem{KK2} Katz, K.; Katz, M.: When is .999\ldots{} less than 1?
\emph{The Montana Mathematics Enthusiast} \textbf{7} (2010), No.~1,
3--30.

See http://arxiv.org/abs/arXiv:1007.3018


\bibitem{KK11b} Katz, K.; Katz, M.: Cauchy's continuum.
\emph{Perspectives on Science} \textbf{19} (2011), no.~4, 426-452.
See http://arxiv.org/abs/1108.4201 and

http://www.mitpressjournals.org/doi/abs/10.1162/POSC\_a\_00047


\bibitem{KK11d} Katz, K.; Katz, M.: Meaning in classical mathematics:
is it at odds with Intuitionism?  \emph{Intellectica} \textbf{56}
(2011), no.~2, 223--302.  See

http://arxiv.org/abs/1110.5456



\bibitem{KK11a} Katz, K.; Katz, M.: A Burgessian critique of
nominalistic tendencies in contemporary mathematics and its
historiography.  \emph{Foundations of Science} \textbf{17} (2012),
no.~1, 51--89.  See http://dx.doi.org/10.1007/s10699-011-9223-1

and {http://arxiv.org/abs/1104.0375}




\bibitem{KK11c} Katz, K.; Katz, M.: Stevin numbers and reality.
\emph{Foundations of Science} \textbf{17} (2012), no.~2, 109-123.  See
http://dx.doi.org/10.1007/s10699-011-9228-9

and http://arxiv.org/abs/1107.3688



\bibitem{KL} Katz, M.; Leichtnam, E.: Commuting and non-commuting
infinitesimals.  \emph{American Mathematical Monthly} (to appear).


\bibitem{KSS} Katz, M.; Schaps, D.; Shnider, S.: Almost Equal: The
Method of Adequality from Diophantus to Fermat and Beyond,
\emph{Perspectives on Science} \textbf{20} (2012), to appear.


\bibitem{KS1} Katz, M.; Sherry, D.: Leibniz's infinitesimals: Their
fictionality, their modern implementations, and their foes from
Berkeley to Russell and beyond.  \emph{Erkenntnis} (online first), see
http://dx.doi.org/10.1007/s10670-012-9370-y

and http://arxiv.org/abs/1205.0174

\bibitem{KS2} Katz, M.; Sherry, D.: Leibniz's laws of continuity and
homogeneity.  \emph{Notices of the American Mathematical Society}, to
appear.


\bibitem{KT} Katz, M.; Tall, D.: The tension between intuitive
infinitesimals and formal mathematical analysis.  Chapter in: Bharath
Sriraman, Editor.  Crossroads in the History of Mathematics and
Mathematics Education.  \emph{The Montana Mathematics Enthusiast
Monographs in Mathematics Education} \textbf{12}, Information Age
Publishing, Inc., Charlotte, NC, 2011, pp.~71-89.  See

http://arxiv.org/abs/1110.5747 and

http://www.infoagepub.com/products/Crossroads-in-the-History-of-Mathematics




\bibitem{Ke} Keisler, H. Jerome: Elementary Calculus: An Infinitesimal
Approach.  Second Edition.  Prindle, Weber \& Schimidt, Boston, 1986.


\bibitem{Ke94} Keisler, H. Jerome: The hyperreal line.  Real numbers,
generalizations of the reals, and theories of continua, 207--237,
Synthese Lib., 242, Kluwer Acad. Publ., Dordrecht, 1994.


\bibitem{Kl08} Klein, F.: Elementary Mathematics from an Advanced
Standpoint.  Vol. I.  Arithmetic, Algebra, Analysis.  Translation by
E. R. Hedrick and C. A. Noble [Macmillan, New York, 1932] from the
third German edition [Springer, Berlin, 1924].  Originally published
as Elementarmathematik vom h\"oheren Standpunkte aus (Leipzig, 1908).



\bibitem{Lau87} Laugwitz, D.: Infinitely small quantities in Cauchy's
textbooks.  \emph{Historia Mathematica} \textbf{14} (1987), no.~3,
258--274.


\bibitem{Lau89} Laugwitz, D.: Definite values of infinite sums:
aspects of the foundations of infinitesimal analysis around 1820.
\emph{Archive for History of Exact Sciences} \textbf{39} (1989),
no.~3, 195--245.


\bibitem{Lau92} Laugwitz, D.: Early delta functions and the use of
infinitesimals in research.  \emph{Revue d'histoire des sciences}
\textbf{45} (1992), no.~1, 115--128.



\bibitem{Li} Lightstone, A. H.: Infinitesimals.  \emph{American
Mathematical Monthly} \textbf{79} (1972), 242--251.


\bibitem{LW} Lightstone, A. H.; Wong, K.: Dirac delta functions via
nonstandard analysis.  \emph{Canadian Mathematical Bulletin}
\textbf{18} (1975), no.~5, 759--762.


\bibitem{Lo} {\L}o{\'s}, J.: Quelques remarques, th\'eor\`emes et
probl\`emes sur les classes d\'efi\-nissables d'alg\`ebres.  In
Mathematical interpretation of formal systems, {98--113},
North-Holland Publishing Co., Amsterdam, 1955.
%



\bibitem{Lut82} L\"utzen, J.: The prehistory of the theory of
distributions.  \emph{Studies in the History of Mathematics and
Physical Sciences} \textbf{7}.  Springer-Verlag, New York-Berlin,
1982.
%
%

\bibitem{Lu62} Luxemburg, W.: Nonstandard analysis.  Lectures on
A. Robinson's Theory of infinitesimals and infinitely large numbers.
Pasadena: Mathematics Department, California Institute of Technology'
second corrected ed., 1964.



\bibitem{Me} M\'eray, H. C. R.: Remarques sur la nature des
quantit\'es d\'efinies par la condition de servir de limites \`a des
variables donn\'ees, \emph{Revue des soci\'eti\'es savantes des
d\'epartments, Section sciences math\'ematiques, physiques et
naturelles, 4th ser.}, \textbf{10} (1869), 280--289.


\bibitem{Ne} Nelson, E.: Internal set theory: a new approach to
nonstandard analysis.  \emph{Bulletin of the American Mathematical
Society} \textbf{83} (1977), no.~6, 1165--1198.


\bibitem{No2} Novikov, S. P.: The second half of the 20th century and
its conclusion: crisis in the physics and mathematics community in
Russia and in the West.  Amer. Math. Soc. Transl. Ser. 2, 212,
Geometry, topology, and mathematical physics, 1--24, Amer. Math. Soc.,
Providence, RI, 2004.  (Translated from Istor.-Mat. Issled. (2)
No. 7(42) (2002), 326--356, 369; by A. Sossinsky.)

\bibitem{Pei} Peirce, C. S.: Three Grades of Clearness.  In ``The
Logic of Relatives'', \emph{The Monist} \textbf{7} (1897),
pp. 161--217.

\bibitem{Per} P\'eraire, Y.: A mathematical framework for Dirac's
calculus.  \emph{Bulletin of the Belgian Mathematical Society. Simon
Stevin} \textbf{13} (2006), no.~5, 1007--1031.


\bibitem{Ro61} Robinson, A.: Non-standard analysis.
\emph{Nederl. Akad. Wetensch. Proc. Ser. A} \textbf{64} =
\emph{Indag. Math.} \textbf{23} (1961), 432--440 [reprinted in
Selected Works, see item~\cite{Ro79}, pp.~3-11]



\bibitem{Ro67} Robinson, A.: The metaphysics of the calculus.  In
Problems in the philosophy of mathematics, Imre Lakatos, Ed.
North-Holland Publ. Co., Amsterdam, 1967, pp.~28-46.




\bibitem{Ro79} Robinson, A.: Selected papers of Abraham
Robinson. Vol. II.  Nonstandard analysis and philosophy. Edited and
with introductions by W. A. J. Luxemburg and S. Körner. Yale
University Press, New Haven, Conn., 1979.


\bibitem{Roq} Roquette, P.: Numbers and models, standard and
nonstandard.  \emph{Math Semesterber} \textbf{57} (2010), 185--199.


\bibitem{SL} Schmieden, C.; Laugwitz, D.: Eine Erweiterung der
Infinitesimalrechnung.  \emph{Mathematische Zeitschrift} \textbf{69}
(1958), 1--39.


\bibitem{She87} Sherry, D.: The wake of Berkeley's Analyst:
\emph{rigor mathematicae}? \emph{Studies in History and Philosophy of
Science} \textbf{18} (1987), no.~4, 455--480.

\bibitem{Si} Sinaceur, H.: Cauchy et Bolzano.  \emph{Revue d'Histoire
des Sciences et de leurs Applications} \textbf{26} (1973), no.~2,
97--112.


\bibitem{Sk} Skolem, Th.: \"Uber die Nicht-charakterisierbarkeit der
Zahlenreihe mittels endlich oder abz\"ahlbar unendlich vieler Aussagen
mit ausschliesslich Zahlenvariablen.  \emph{Fundamenta Mathematicae}
\textbf{23}, 150-161 (1934).


\bibitem{Sm} Smithies, F.: Review of L\"utzen (item \cite{Lut82}
above), 1984.

See http://www.ams.org/mathscinet-getitem?mr=667854


\bibitem{Str} Stroyan, K.: Uniform continuity and rates of growth of
meromorphic functions. Contributions to non-standard analysis
(Sympos., Oberwolfach, 1970), pp. 47--64.  \emph{Studies in Logic and
Foundations of Math.}, Vol. 69, North-Holland, Amsterdam, 1972.

\bibitem{Ta} Takeuti, G.: Dirac space. \emph{Proceedings of the Japan
Academy} \textbf{38} (1962), 414--418.


\bibitem{To} Todorov, T. D.: A nonstandard delta function.
\emph{Proceedings of the American Mathematical Society} \textbf{110}
(1990), no. 4, 1143--1144.

\bibitem{Ya} Yamashita, H.: Comment on: ``Pointwise analysis of scalar
fields: a nonstandard approach'' [J. Math. Phys. 47 (2006), no.~9,
092301; 16 pp.].  \emph{Journal of Mathematical Physics} \textbf{48}
(2007), no.~8, 084101, 1 page.


\end{thebibliography}
\end{document}